\documentclass[11pt]{article} 
 
\usepackage{latexsym} 
\usepackage{amsfonts} 
\usepackage{amsmath} 
\usepackage{amsthm} 
\usepackage{mathrsfs} 
\usepackage{dsfont} 
\usepackage{bbold} 
\usepackage[english]{babel} 
\usepackage{caption} 
\usepackage{epsfig} 
\usepackage{float} 
\usepackage{subfigure} 
\usepackage{psfrag} 
\usepackage{graphicx} 
\usepackage{epsfig} 
\usepackage{tensor}
\usepackage{hyperref}
\usepackage{amssymb}

\DeclareMathOperator{\meas}{meas}
\DeclareMathOperator{\diag}{diag} 
\DeclareMathOperator{\dist}{dist} 
\DeclareMathAlphabet{\mathpzc}{OT1}{pzc}{m}{it}
 
\newtheorem{theorem}{Theorem} 
\newtheorem*{prop*}{Theorem}

\newtheorem{defi}[theorem]{Definition} 
\newtheorem{lemma}[theorem]{Lemma} 
 
\newtheorem{rmk}[theorem]{Remark}

\numberwithin{equation}{section}
\numberwithin{theorem}{section}

\newcommand{\ZZZ}{\mathds{Z}} 
\newcommand{\CCC}{\mathds{C}} 
\newcommand{\NNN}{\mathds{N}} 
 
\newcommand{\RRR}{\mathds{R}} 
\newcommand{\TTT}{\mathds{T}} 
\newcommand{\uno}{\mathds{1}} 
 \newcommand{\id}{\mathds{1}} 

\newcommand{\CCCC}{{\mathcal C}} 
\newcommand{\DD}{{\mathcal D}}

\newcommand{\calG}{{\mathtt M}} 
\newcommand{\calH}{{\mathcal H}} 
\newcommand{\calI}{{\mathcal I}} 
\newcommand{\calL}{{\mathcal L}}

\newcommand{\LL}{{\mathcal L}} 
 
\newcommand{\NN}{{\mathcal N}}

\newcommand{\calS}{{\mathcal S}} 
\newcommand{\Gacca}{{\mathtt G}} 
\newcommand{\SSSSGa}{\Lambda_{+}(\Gacca)}
\newcommand{\SSSSG}{\Lambda_{+}}

 \newcommand{\pD}{{\mathtt D}} 
  \newcommand{\pDL}{{\mathtt L}}

\newcommand{\calmR}{{\mathscr R}}

\newcommand{\gota}{{\mathfrak a}}

\newcommand{\gotd}{{\mathfrak d}}

\newcommand{\gotn}{{\mathfrak n}}

\newcommand{\gotA}{\{-1,1\}}

\newcommand{\gotF}{{\mathfrak F}} 
 
\newcommand{\gotH}{{\mathfrak H}}

\newcommand{\gotK}{{\mathfrak K}}

\newcommand{\ol}{\overline} 
 
\newcommand{\Fullbox}{{\rule{2.0mm}{2.0mm}}} 
 
\newcommand{\EP}{\hfill\Fullbox\vspace{0.2cm}} 
 
\newcommand{\io}{\infty} 
\newcommand{\e}{\varepsilon} 
 
\newcommand{\de}{\delta} 
\newcommand{\be}{\beta}

\newcommand{\g}{\gamma} 
\newcommand{\om}{\omega}

\newcommand{\la}{\lambda} 
\newcommand{\f}{\varphi} 
\newcommand{\s}{\sigma} 
 
\newcommand{\del}{\partial}

\newcommand{\av}[1]{\langle #1 \rangle}

\newcommand{\ff}{\boldsymbol{f}}

\newcommand{\ii}{{\rm i}}

\newcommand{\se}{\mathfrak a}

\def\tilde#1{\widetilde{#1}}

\def\ins#1#2#3{\vbox to0pt{\kern-#2 \hbox{\kern#1 #3}\vss}\nointerlineskip} 
 
\begin{document}

\title{\bf A KAM result on compact Lie groups} 
 
\author{\bf Livia Corsi$^{1}$, Emanuele Haus$^{2}$, Michela Procesi$^{3}$
\vspace{2mm} 
\\ \small
$^{1}$Department of Mathematics and Statistics, McMaster University, Hamilton, ON, L8S 4K1, Canada
\\ \small 
$^2$Dipartimento di Matematica, Universit\`a di
Napoli ``Federico II'', Napoli, I-80126, Italy
\\ \small
$^2$Dipartimento di Matematica, Sapienza - Universit\`a di Roma, Roma, I-00185, Italy
\\ \small 
E-mail:  lcorsi@math.mcmaster.ca, emanuele.haus@unina.it, mprocesi@mat.uniroma1.it}
 
\date{} 
 
\maketitle 
 
\begin{abstract}
We describe some recent results on existence of quasi-periodic solutions of Hamiltonian PDEs on compact manifolds. We prove a linear stability result for the non-linear Schr\"odinger equation in the case of $SU(2)$ and $SO(3)$.

\smallskip

\smallskip


\smallskip

\end{abstract} 

  \tableofcontents

\section{Introduction} 
\label{sec.intro} 

In the past forty years there has been a lot of progress in the study of 
many non-linear PDEs which model the propagation of waves.
In this class of equations we can mention for instance the non-linear wave (NLW) equation, the Euler
equations of hydrodynamics and various models deriving from it such as the non-linear Schr\"odinger (NLS)
the Korteweg-de Vries (KdV), the Camassa-Holm equations and many others.

A particularly fruitful research line has been the so-called
``dynamical systems approach'' i.e. the generalisation to infinite dimensional setting
of many ideas and techniques borrowed from the theory of dynamical systems; the key idea
is to look for invariant manifolds on which the dynamics is particularly simple and then
try to obtain some stability result in order to deduce some properties for typical initial data
on the whole phase space.

The behaviour of the solutions is expected to depend strongly on the set in which the ``space variable''
lives; in this paper we will concentrate on the case of a compact Riemannian manifold ${\mathtt M}$, where
one expects a ``recurrent dynamics'' and complicated coexistence of regular and chaotic phenomena.
In particular we shall focus on the problem of existence and stability of quasi-periodic solutions and,
as an example, we will study a forced NLS equation
\begin{equation}\label{nls}
\ii u_{t}-\Delta u+{\mathtt m}u = \e {\mathtt f}(\om t,x,u),
\qquad x\in {\mathtt M} \, .
\end{equation}
Here and henceforth $ \Delta $ denotes the Laplace-Beltrami operator,   ${\mathtt m}>0 $ is the ``mass", 
the parameter $\e >  0 $ is  small, and  the frequency vector is 
$\om\in\RRR^{d}$. 
Concerning regularity we assume that
${\mathtt f}(\varphi, x, u ) \in C^q(\TTT^d\times{\mathtt M}\times \CCC ;\CCC)$ in the real sense
(namely as a function of ${\rm Re}(u), {\rm Im}(u)$),  for some $ q $ large enough.

We will describe an existence result of quasi-periodic solutions for equation \eqref{nls} above when $\mathtt M$ is a homogeneous space w.r.t. a compact Lie group (see Theorem \ref{teoremone}). Then we shall restrict our attention to the case of $SU(2)$ or $SO(3)$ and prove a linear stability result (see Theorem \ref{thm.reduci}). In Section \ref{salaminchia} we discuss the extension to spherical varieties of rank one.

\medskip

\noindent
\emph{The Newton scheme.}
Passing to the Fourier representation for the space variables, namely
$$
u(x,t)=\sum_j u_j(t)\phi_j(x)
$$ 
where $j$ runs in a countable index set and the $\phi_j$'s are the eigenfunctions of $\Delta$,
\eqref{nls} can be seen as an infinite dimensional forced dynamical system which has
an elliptic fixed point at $\e=0$. A very natural question is whether there are solutions
which synchronise with the forcing, provided that the forcing frequency $\om$ is sufficiently
non-resonant w.r.t. the linear frequencies, i.e. the eigenvalues of the operator $\Delta$.

Although \eqref{nls} is a simplified problem with respect to the autonomous case,
it still contains some of the main difficulties
that one has to deal with and a full understanding of even this simplified case is an open problem.
Indeed, as a first na\"\i ve attempt, one reduces the search for quasi-periodic solutions to \eqref{nls}
to an implicit function problem
$$
F(u):=F(u,\om)=0
$$
and may try to solve it by perturbation theory; however the
linearised operator at $\e=0$ is $\ii\del_t-\Delta+{\mathtt m}$ and its inverse is unbounded so
that one cannot apply the Implicit Function Theorem: this is known as the ``small divisor problem''.
 
\begin{figure}[ht] 
\centering 
\ins{225pt}{-20pt}{$F(u)$}
\ins{212pt}{-143pt}{$u_0$}
\ins{172pt}{-143pt}{$u_1$}
\ins{149pt}{-143pt}{$u_2$}
\includegraphics[width=3in]{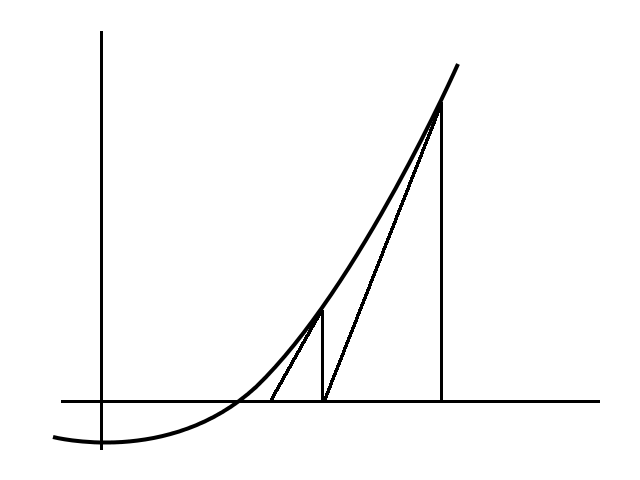} 
\vskip.2truecm 
\caption{Three steps of the Newton algorithm
$u_{n+1} := u_n - (F'(u_n,\om))^{-1} [F(u_n)] $ }
\label{fig.newton} 
\end{figure}

In order to handle this, one typically uses a recursive Newton-like scheme, which is based on the
invertibility of the equation linearised at a sequence of {\em approximate solutions} $u_n$,  
 see Figure \ref{fig.newton}.


This in turn can be seen as a ``non-resonance'' condition on the frequency $\om$: 
indeed one can produce an abstract ``Nash-Moser'' scheme (see for instance \cite{BBP,BCP})
which says that if $\om$ is such that at each step $n$ of the scheme the operator $(F'(u_n,\om))^{-1}$ is well-defined
and bounded from $H^{s+\mu}$ to $H^{s}$ for some $\mu$, then a solution of \eqref{nls} exists.
Then the problem reduces to proving that such set of parameters $\om$ is non-empty, or even better
that it has asymptotically full measure.

If we impose some symmetry such as a Hamiltonian structure, the linearised operator
$F'(u,\om)$ is self-adjoint and it is easy to obtain some information on its eigenvalues,  
implying its invertibility with bounds
on the $L^2$-norm of the inverse for ``most" parameters $ \om $. However  this information
is not enough to prove the convergence of the algorithm: one needs
estimates on the high Sobolev norm of the inverse, which  do not follow only
from bounds on the eigenvalues. 

Naturally, if $F'(u,\om)$ were diagonal, passing from $L^2$ to $H^s$ norm would be trivial,
but the problem is that the operator which diagonalises $F'(u,\om)$ may not be bounded in
$H^s$. The property of an operator to be diagonalisable via a ``smooth'' change of variables
is known as \emph{reducibility} and in general is connected to the fact that the matrix is regular
semi-simple, namely its eigenvalues are distinct.
When dealing with infinite dimensional matrices, one also has to give quantitative estimates
on the difference between two eigenvalues: this is usually referred to as the \emph{second order
Mel'nikov condition}, since it can be seen as a condition on $\om$. However in general this condition
cannot be imposed because the eigenvalues of $\Delta$ are multiple and actually have
unbounded multiplicity. Naturally one does not need to diagonalise a matrix in order to invert it,
and indeed there are various existence results which have been proved in the case of multiple
eigenvalues; however, this tends to be technically quite complicated and needs a deep understanding
of the harmonic analysis on the manifold ${\mathtt M}$.

\medskip
 
\noindent
\emph{Some literature.}
The first existence (and stability) results dealt with autonomous Hamiltonian PDEs and were obtained
by Kuksin \cite{K1}, P\"oschel \cite{KP,Po2}, Wayne
\cite{W1} who studied the NLS and NLW equations on the interval $[0,\pi]$ where the eigenvalues of $\Delta$ are
simple and one can easily impose the second order Mel'nikov conditions. Thanks to this
diagonalisation procedure, they were able to obtain some information on the linear stability: in particular
they showed that the eigenvalues are purely imaginary. Their approach was an infinite dimensional
generalisation of the classical KAM algorithm for elliptic tori (see for instance \cite{M,Po3}).

Since these results dealt with autonomous equations, another problem was the so-called \emph{frequency
modulation}, namely the fact that there are no external parameters on which to impose the non-resonance
conditions and one needs to ``extract them'' from the nonlinearity, in general by means of Birkhoff normal form.

Later on, these KAM techniques were further generalised by Chierchia-You to the case of NLW with multiple
eigenvalues but with bounded multiplicity, for instance when $\mathtt M$ is the unit circle $\TTT$.

A more direct approach was proposed by Craig-Wayne \cite{CW}, who dealt with an analytic setting;
they used a Lyapunov-Schmidt decomposition
in order to ``extract the parameters'' and a Newton scheme to solve the small divisor problem.
In order to get the needed estimates on
$(F'(u))^{-1}$ (also called \emph{Green function estimates} by analogy with the Anderson localisation problem)
they developed a technique inspired by the methods of Fr\"olich-Spencer \cite{FS}.
However their result was limited to the case of periodic solutions.

This set of techniques was extended to the case of quasi-periodic solutions by Bourgain
\cite{B1,B3,B5} who was able to deal also with the case ${\mathtt M}=\TTT^d$.
Then Berti-Bolle \cite{BB1,BB2} were able to generalise Bourgain's techniques to the case of Sobolev regularity,
considering also a multiplicative potential (the previous results dealt with a simplified model where the potential
is non-local).

The reason why these results are confined to tori 
 is that their proofs  require specific properties  of the eigenvalues, while the eigenfunctions 
must be the exponentials or, at least,  strongly ``localised close to exponentials".
In the  paper \cite{BP} Berti-Procesi proved existence of periodic solutions for NLW and NLS on any compact  Lie group
or manifold homogenous  with respect to a compact Lie group
 and finally Berti-Corsi-Procesi \cite{BCP} extended this result to the case of quasi-periodic solutions.

Note that all these results obtained via Newton method
do not give any information on the linear stability of the solution, which is a completely non-trivial problem
since the second order Mel'nikov condition is obviously violated already on $\TTT^d$ with $d\geq2$.

The first reducibility results on $\TTT^d$ with $d\geq 2$ are due to Eliasson-Kuksin \cite{EK1,EK} who were able to prove
 linear stability of the quasi-periodic solutions of NLS. The main ingredients of their proofs are the following:
first they reduce to a time-independent block diagonal matrix and then they impose the second
order Mel'nikov condition
between the eigenvalues of different blocks.
In order to show that the set of parameters has positive measure they need to study carefully the
asymptotics of the eigenvalues (the so called T\"oplitz-Lipschitz condition). 
We mention also the papers \cite{GXY,PX,PP} which make use of the conservation of momentum
in order to fully diagonalise the matrix.

Very recently a combination of the two approaches has been developed by Baldi-Berti-Montalto
\cite{BBM1,BBM3,BBM2} in order to prove existence and stability for fully non-linear perturbations of
the KdV equation; see also \cite{FP} for the case of the NLS equation.
We believe that this latter approach may be very fruitful since it decouples completely the
existence and reducibility problems; note however that the strategy used so far in order to deal with
unbounded perturbations works only in one-dimensional cases.

In this paper, for $\mathtt M=SU(2), SO(3)$, we shall prove by means of a KAM reducibility scheme the linear stability of the quasi-periodic solutions whose existence has been proved in \cite{BCP} in a more general and abstract setting. 
 At a formal level a KAM reducibility scheme starts with a matrix $\mathcal L_\e$ of the form $D+ \e T$   where $D$ is diagonal with distinct eigenvalues  and $T$ is bounded in some appropriate norm. Then one step of the  scheme provides a change of variables which conjugates $\mathcal L_\e$ to $D_1+\e^2 T_1$ where again $D_1$ is  diagonal and $T_1$ bounded. Iterating this procedure one diagonalises the matrix.
In our case we  diagonalise the operator linearised along the solution via a smooth, time quasi-periodic change of variables on the phase space; then we obtain the linear stability by explicitly checking that the eigenvalues of such linearised operator are purely imaginary.

\subsection{Main results}

Let us consider \eqref{nls} where ${\mathtt M}$ is a compact Lie group or manifold
which is homogeneous w.r.t. a compact Lie group (namely there exists a compact Lie
group which acts on ${\mathtt M}$ transitively and differentiably). Assume that \eqref{nls} is Hamiltonian, i.e.
\begin{equation}\label{HamNLS}
{\mathtt f}(\om t,x,u) = \partial_{\ol{u}} H(\om t,x,u) \, , 
\quad H( \varphi,x,u) \in \RRR \, , \ \forall u \in \CCC
\end{equation}
 with real Hamiltonian  
\begin{equation}\label{HamNLS2}
H(\om t,x,u) = \overline{H(\om t,x,u)}\ .
\end{equation}

We assume that the frequency $\om$ has a fixed Diophantine direction,  namely 
\begin{equation}\label{omegabar}
\om=\la\tilde\om,\qquad \la\in\mathcal I:=[1/2,3/2],\qquad
|\tilde\om |_1:= {\mathop\sum}_{p=1}^d|\tilde\om_p| \le1,
\end{equation}
for some fixed Diophantine vector $\tilde\om$, i.e. which
satisfies 
\begin{equation}\label{dioph}
|\tilde\om\cdot l|\ge 2\g_{0} |l|^{-d}, \quad \forall\,l\in\ZZZ^{d}\setminus\{0\},
\end{equation}
for some positive $\g_0$.
The search for quasi-periodic solutions of \eqref{nls} reduces to finding solutions $ u(\varphi, x) $  of 
\begin{equation}\label{NLSNLW}
\ii \om \cdot \partial_\varphi u -\Delta u+{\mathtt m}u = \e{\mathtt f}(\varphi ,x,u) \, ,  
\end{equation}
in some Sobolev space $ H^s $ of both variables $ ( \varphi , x ) $.

It is convenient to ``double'' the NLS equation \eqref{NLSNLW}, namely consider the 
\emph{vector} NLS operator
\begin{equation}\label{vnls}
F(\e, \la, u^+,u^-):=\left\{
\begin{aligned}
\ii \la \tilde\om\cdot\del_{\f}u^{+}-\Delta u^{+}+{\mathtt m}u^+-\e \gotF(\f,x,u^{+},u^{-})\\
-\ii \la \tilde\om\cdot\del_{\f}u^{-}-\Delta u^{-}+{\mathtt m}u^--\e {\gotH}(\f,x,u^{+},u^{-})
\end{aligned}
\right.
\end{equation}
on the space $ H^s(\TTT^d\times \calG)\times H^s(\TTT^d\times \calG)$,
where $\gotF(u,v),\gotH(u,v) $ are two extensions of class $ C^q (\TTT^d \times {\mathtt M \times \CCC^2; \CCC}) $
(in the real sense)  of ${\mathtt f}(u) $
such that $\gotF(u,\ol{u})=\ol{\gotH(u,\ol{u})}={\mathtt f}(u)$ and
$\del_u \gotF(u,\ol{u})=\del_v \gotH(u,\ol{u})\in\RRR$,
$\del_{\ol u} \gotF(u,\ol{u})=\del_{\ol u}\gotH(u,\ol{u})=
\del_{\ol v} \gotF(u,\ol{u})=\del_{\ol v}\gotH(u,\ol{u})=0$ and
$\del_v \gotF(u,\ol{u})=\ol{\del_u\gotH(u,\ol{u})}$; see for instance \cite{BB1,BCP}.

Note that \eqref{vnls} reduces to \eqref{NLSNLW} on the invariant subspace
$$
{\cal U} := \{u=(u^{+},u^{-})\in H^{s}\times H^{s}\;:\; u^{-}=\ol{u^{+}}\}.
$$

The following result has been proved in \cite{BCP}.

\begin{theorem}[{\bf Existence}]\label{teoremone}
Let  ${\mathtt M}$ be a compact Lie group or a manifold homogeneous w.r.t. a compact Lie group,
consider the vector NLS equation $F(\e,\la,u^+,u^-)=0$ where $F$ is the non-linear operator in \eqref{vnls}
and assume   \eqref{omegabar}-\eqref{dioph}. Then there are some large numbers $s_1,  q ,S \in \RRR $ 
such that, for any $ {\mathtt f} \in C^{q}$ and for all $\e\in[0,\e_{0})$ with $\e_{0}>0$ small enough,
there is a map
$$
u_\e\in C^{1}(\calI,H^{s_1}),\qquad
\sup_{\la\in\calI}\|u_\e (\la)\|_{s_1}\to0,\mbox{ as }\e\to0,
$$
and a set $\CCCC_{\e}\subseteq\calI$, satisfying $ {\meas(\CCCC_{\e})}=1-O(\e^{1/S})$, 
such that, for any $\la\in\CCCC_{\e}$, $w_\e(\la):=(u_\e (\la),\ol{u}_\e(\la))$ is a solution of
\eqref{vnls}, with $\om=\la\tilde\om$. Moreover if 
 $ {\mathtt f} \in C^{\io}$ then
$u_\e (\la) $ is of class $ C^{\io} $ both in time and space. Finally if $\mathtt f$ is central on $\calG$,
i.e.
\begin{equation}\label{merdecentrali}
\mathtt f(\om t,x,u)=\mathtt f(\om t, g^{-1}x g,u)\,,\quad \forall g\in \calG
\end{equation}
then $u_\e(\la)$ is central.
\end{theorem}

Actually the last sentence is not explicitly stated in \cite{BCP} but it follows directly from
\cite{BCP}-Corollary 2.17.

The proof of Theorem \ref{teoremone} relies on an abstract Nash-Moser scheme on sequence
spaces; as explained above,
the convergence of such scheme only requires ``tame'' 
estimates of the inverse in high Sobolev norm. Following \cite{BB1,BB2}, such estimates
have been obtained by means of a {\it multiscale analysis}.
Roughly speaking, it is a way to prove an off-diagonal decay (see Definition \ref{def:Ms})
for the inverse of a finite-dimensional invertible matrix with off-diagonal decay, by using information on the
invertibility (in high norm) of a sufficient number of principal minors of order $N$ much
smaller than the dimension of the matrix.
In applying these ideas to the case of Lie groups, two key points concern
\begin{enumerate}
\item  the matrix representation of a multiplication operator
$ u \mapsto b u $,
\item the properties of the eigenvalues of the Laplace-Beltrami operator. 
\end{enumerate}
The multiplication rules for the eigenfunctions imply that 
the operator of multiplication by a Sobolev function $ b \in H^s ({\mathtt M})$ is
represented -- in the eigenfunction basis -- as a block matrix with off-diagonal
decay, as  stated precisely in Lemma \ref{moltiplicazione.matrici} (proved in \cite{BP}). 
 The block structure of this matrix takes into account the (large) multiplicity of the
degenerate eigenvalues of $ \Delta $ on ${\mathtt M} $. 
This in principle could be a problem because one cannot
hope to achieve any off-diagonal decay property for the matrices
 restricted to such blocks; actually, we can only control the $L^2$-operator norm on these blocks, but this is enough to prove the existence result.

Concerning item 2, the eigenvalues of the Laplace-Beltrami operator on a Lie group
are very  similar to those on a torus. This enables one to prove 
``separation properties'' of clusters of singular/bad sites (i.e. Fourier indices corresponding to a small eigenvalue) 
\textit{\`a la} Bourgain \cite{B3,B5}. 
Thanks to the off-diagonal decay property discussed in item 1, such ``resonant'' clusters 
 interact only weakly. 

Under the hypotheses of Theorem \ref{teoremone}, set
\begin{equation}\label{line}
\LL_\e:=
(-\Delta+{\mathtt m})\sigma_3 -\e T
\end{equation}
where
\begin{equation}\label{Ti}
T=T(w_\e):=
\begin{pmatrix}
D_{u^+}{\gotF}(\om t,x,u_\e(\la),\ol{u}_\e(\la)) & -D_{u^-}{\gotF}(\om t,x,u_\e(\la),\ol{u}_\e(\la)) \cr
D_{u^+}{\gotH}(\om t,x,u_\e(\la),\ol{u}_\e(\la)) & -D_{u^-}{\gotH}(\om t,x,u_\e(\la),\ol{u}_\e(\la))
\end{pmatrix}
\end{equation}
and $\s_3$ is the third Pauli matrix, namely
\begin{equation}\label{pauli}
\sigma_3=\begin{pmatrix}
\uno & 0 \cr
0 & -\uno
\end{pmatrix},
\end{equation}
i.e. $-\ii\LL_\e$ is the vector field linearised at the solution. In the present paper we shall prove the following result.

\begin{theorem}[{\bf Linear Stability}]\label{thm.reduci}
Assume that $\calG=SU(2),SO(3)$ and $\mathtt f$ is central on $\calG$ (see \eqref{merdecentrali}). Then
under the same assumptions of Theorem \ref{teoremone}, possibly with smaller $\e_0$ and larger $q$,  there exist $\alpha>0$, $s_2\leq s_1-\alpha$  and a subset
$\calS\subseteq\CCCC_\e$ such that for $\la\in\calS$ the equation \eqref{vnls} linearised at
the solution $w_\e(\la)$ is reducible with a change of variables in $H^{s_2}\times H^{s_2}$. More precisely 
\begin{equation}\label{pillola}
{\rm meas}(\calS)\to1 \quad \text{as}\quad \e\to0
\end{equation}
and for all $\lambda\in\calS$ there
exists a  quasi-periodic close-to-identity
change of variables $h=\Psi(\om t) v$ which reduces the linearised vector NLS equation
\begin{equation}\label{lin.nls}
 h_t + \ii \LL_\e h =0
\end{equation}
to
\begin{equation}\label{ridotto}
 v_t + \ii \DD v=0
\end{equation}
with $\DD$ a diagonal and time-independent linear operator whose eigenvalues are explicitly given in formula
\eqref{eq:4.4}. Finally for all $\varphi\in \TTT^d$, one has
\begin{equation}\label{poropo} \| [\Psi(\varphi)]^{-1}h-h\|_{s_2}\leq C \e^a (1 +\|u_\e(\la)\|_{s_2+\alpha})\|h\|_{s_2}\,,\quad \forall h \in H^{s_2}(\calG,\CCC)  
\end{equation}
for some $a\in (0,1)$ and some $\varphi$-independent constant $C$.
Finally one has
\begin{equation}\label{scaramouche}
1- K\e^a (1 +\|u_\e(\la)\|_{s_2+\alpha}) \le \frac{\|h(t)\|_{s_2}}{\|h(0)\|_{s_2}}\le 1+ K\e^a (1 +\|u_\e(\la)\|_{s_2+\alpha}),
\end{equation}
for some  constant $K$.
\end{theorem}

We confine ourselves to the case of $SU(2),SO(3)$ in order
to have a precise control on the differences of the eigenvalues $\mu_j$ of $-\Delta$; see Section \ref{sec:harmonic}.
This in turn will allow us to impose the second order Mel'nikov conditions. Note that, differently from the existence result, here we restrict ourselves to central functions in order to avoid having to deal with multiple eigenvalues. In principle, one could weaken this restriction and obtain a block diagonal, time-independent matrix $\mathcal D_\e$. However at the moment we are not able to prove the convergence of the resulting KAM scheme and actually it is not even clear to us whether this is a technical or a substantial problem.


\medskip

\section{The functional setting}

A compact manifold  $ \mathtt M $ which is homogeneous w.r.t. a compact Lie
group is, up to an isomorphism,  diffeomorphic to 
\begin{equation}\label{MGN}
\mathtt M =  G/ N \, , \quad G := \Gacca\times\TTT^{\gotn} \, , 
\end{equation}
where $ \Gacca $ is a simply connected compact Lie group, $ \TTT^{\gotn} $ is a torus  and $ N $ is a closed 
subgroup of $  G  $. Then a function on $ \mathtt M  $ can be seen 
as a function defined on $ G $ which is invariant under the action of $ N $, and  
the space $H^s(\mathtt M,\CCC)$ can be identified 
with the subspace 
\begin{equation}\label{subsp}
\widehat H^s := 
\widehat H^s(G ,\CCC):= 
\Big\{ u \in H^s(G)\,:\; u(x)= u(x g )\,,\ \ \forall  x \in G = \Gacca\times\TTT^{\gotn}, g \in N \Big\}. 
\end{equation}
Moreover, the Laplace-Beltrami operator on  $ \mathtt M $ can be identified with the Laplace-Beltrami operator 
on the Lie group $ G  $, acting on functions invariant under $ N $ (see Theorem 2.7, \cite{BP}). 
Then we ``lift" the equation \eqref{nls} on  $ G $ and we use harmonic analysis on Lie groups.

\subsection{Analysis on  Lie groups}\label{sec:harmonic}

Any simply connected compact Lie group $\Gacca$ is the product of a finite number of
simply connected 
Lie groups of simple type (which are classified and come in a finite number of families).

Let $\Gacca$ be of simple type with dimension $\gotd$ and rank $r$.
Denote by ${\mathtt w}_{1},\ldots,{\mathtt w}_{r}\in \RRR^r$ the fundamental weights of $\Gacca$
and consider the cone of dominant  weights
$$
\Lambda_+=\SSSSGa:=\Big\{j=\sum_{p=1}^{r}j_{p}{\mathtt w}_{p}\,:\,
j_{p}\in\NNN \Big\} \subset \Lambda := \Big\{j=\sum_{p=1}^{r}j_{p}{\mathtt w}_{p}\,:\,
j_{p}\in \ZZZ \Big\} \, . 
$$
Note that $\SSSSGa$ index-links the finite dimensional irreducible representations of $\Gacca$.

Given an irreducible unitary representation  $(R_{V_{j}},V_{j})$ of $\Gacca$ we denote by
$\ff_{j}(x)$ the (unitary) matrix associated to it, i.e.
\begin{equation}\nonumber
(\ff_{j}(x))_{h,k}=\langle R_{V_{j}}(x)v_{h},v_{k}\rangle,
\qquad v_{h},v_{k} \in V_{j} \, , 
\end{equation}
where $ (v_h)_{h = 1, \ldots, {\rm dim} V_j} $
is an orthonormal basis of the finite dimensional euclidean space $ V_j $ with scalar product
$ \langle \cdot, \cdot \rangle $.
Then the eigenvalues and the eigenfunctions of the
Laplace-Beltrami operator $-\Delta$ on $\Gacca$ are
\begin{equation}\label{autov.lapla}
\mu_{j}:=|j+\rho|_2^{2}-|\rho|_2^{2},
\qquad
\ff_{j,\s}(x),
\quad
x\in \Gacca,
\quad
j\in\SSSSGa,
\quad
\s=1,\ldots,\mathpzc m_{j},
\end{equation}
where 
$\rho:=\sum_{i=1}^{r}{\mathtt w}_{i} $,  $ |\cdot|_2$ denotes the Euclidean norm on $\RRR^{r}$, 
and $\mathpzc m_j=(\dim V_j)^2$ 
satisfies $\mathpzc m_{j}\le |j+\rho|_2^{\gotd-r}$.

Denote by  $\mathcal N_j$  the eigenspace of $ -\Delta $ corresponding to $\mu_j$.
The Peter-Weyl theorem implies the orthogonal decomposition
$$
L^{2} (\Gacca)=\bigoplus_{j\in \SSSSGa} 
\NN_{j} \, . 
$$
If we denote the \emph{central character} of a representation by $\chi_j(x):=\mathrm{tr}(R_{V_j}(x))$, we have that $\{\chi_j\}_{j\in \SSSSGa}$ is a Hilbert basis for the subspace of $L^2(\Gacca)$ formed by the central functions defined in \eqref{merdecentrali}.
\begin{rmk}\label{fankulo}
Note that the multiplicity of an eigenvalue $\mu$ is given by
$$
\sum_{j\,:\,\mu_j=\mu}\mathpzc m_j.
$$
If we reduce to the central functions
we have $\mathpzc m_j=1$; in the case of rank $1$ this implies that the eigenvalues are
simple. 
\end{rmk}

 If $G=SU(2)$ the rank is $1$, the fundamental weight is $\mathtt w_1=(1/4,-1/4)$ and the dominant weights are $j=(m/4,-m/4)$, $m\in\NNN$ so we can identify $\SSSSG$ with $\NNN/\sqrt{8}$. Then the eigenvalues of $-\Delta$ on $SU(2)$ are
\begin{equation}
\left((j+\rho)^2-\rho^2\right)\in\frac{\NNN}{8}
\end{equation}
with $j\in\Lambda_+$ and $\rho=1/\sqrt{8}$. Finally, all the unitary representations $R_{V_j}$ of $SU(2)$ are self-dual (i.e. $\ol{R_{V_j}}=R_{V_j}$), so that the central characters $\chi_j$ are real.

The orthogonal group $SO(3)=SU(2)/\{\pm\uno\}$ is also a homogeneous space and the 
indices of $\Lambda_+(SO(3))$ are half of the indices of $\Lambda_+(SU(2))$. In this case the dominant
weights are $j=(m/2,-m/2)$, $m\in\NNN$ so that $\SSSSG(SO(3))$ is identified with $\NNN/\sqrt{2}$.
From now on we shall consider only $G=SU(2),SO(3)$.

\subsection{Sequence spaces}\label{seq.sp}

The Sobolev space $ H^s(\TTT^d\times G)\times H^s(\TTT^d\times G)$
 can be identified with a sequence space as follows. We start by introducing an index set
\begin{equation}\nonumber 
\mathfrak K := \mathfrak I\times\{-1,1\}:=\ZZZ^d\times \Lambda_+\times \{-1,1\}
\end{equation}
where $\Lambda_+ \subseteq \rho \NNN $ with $\rho=\frac{1}{\sqrt{8}}$.
Given $k\in \mathfrak K$ we denote
\begin{equation}\label{moduliki}
\begin{aligned}
&k=(i,\gota)= (l,j,{\gota}) \in \ZZZ^d\times \rho\ZZZ\times \gotA \\ 
&|k|=|i|:= \max(|l|,|j|),\quad
|l|:=|l|_\io=\max(|l_1|,\ldots,|l_d|)\,.
\end{aligned}
\end{equation}

Finally, for $ k = ( i , \se ), k' = (i' , \se' ) \in{\mathfrak K}$ we denote
\begin{equation}\label{modulo}
\dist(k,k'):=\left\{
\begin{aligned}
&1,\qquad\qquad i=i',\,\se\ne\se',\\
& |i-i'|,\qquad\mbox{otherwise} \, , 
\end{aligned}
\right.
\end{equation}
where $ |i|$ is defined in \eqref{moduliki}.

For $ s \ge 0 $, we define the
(Sobolev) scale of Hilbert sequence spaces
\begin{equation}\nonumber 
H^{s}:=H^{s}({\mathfrak K} ):=\Big\{w=\!\!
\{w_{k}\}_{k\in {\mathfrak K}} \, , \ w_k  \in \CCC  \,:\, 
\|w\|^2_{s}:= 
\sum_{k\in{\mathfrak K}}  |j+\rho| ^{2s}
|w_{k}|^{2}<\infty \Big\}
\end{equation}
%
Similiarly to \eqref{line}--\eqref{Ti} it is convenient to introduce the  following notation: for fixed $i=(l,j)$, $i'=(l',j')$ we set
\begin{equation}\nonumber
M_{\{i\}}^{\{i'\}}:= \{M_{i,\gota}^{i',\gota'}\}_{\gota,\gota'\in \gotA} \,,\quad 
M_{\{i\}}^{\{i'\}}\in {\rm Mat}( 2\times 2,\CCC) \,. 
\end{equation}

\begin{defi}\label{def:Ms} {\bf ($s$-decay norm)}
Fix $s_0>(d+1)/2$. Given a matrix $M$, representing a linear operator on $L^2(\gotK)=H^0(\gotK)$,
we define its $s$-norm   as
\begin{equation}\nonumber
| M|_{s}^{2}:=\sum_{i\in{\mathfrak \ZZZ^{d}\times\rho \ZZZ}}[M(i)]^{2}
\langle i\rangle^{2s}
\end{equation}
where  $\av{i} := \max(1,|i| )$, 
\begin{equation}\nonumber 
[M(i)]:=
\sup_{h-h'=i } \big\|M^{\{h'\}}_{\{h\}} \big\|_0\,.
\end{equation}
If $M=M(\la) $ for $\la\in\calS\subset\RRR$, we define
\begin{equation}\label{2.1}
\begin{aligned}
&|M|_{s,\calS}^{\rm sup}=|M|_{s}^{\rm sup}:=\sup_{\la\in\calS}|M(\la)|_{s}, \quad 
|M|^{\rm lip}_{s,\calS}=|M|^{\rm lip}_{s}:=\sup_{\la_{1}\neq\la_{2}}\frac{|M(\la_{1})-M(\la_{2})|_{s}}{|\la_{1}-\la_{2}|},\\ 
&|M|_{s,\g,\calS}=|M|_{s,\g}:=|M|_{s}^{\rm sup}+\g|M|^{\rm lip}_{s}.
\end{aligned}
\end{equation}
For a Lipschitz family of functions $w(\la)\in H^s(\gotK)$ we define the norm $\|w\|_{s,\gamma}$ exactly in the same way.
Finally, for a Lipschitz function $f:\calI\to \RRR$ we denote by $|f|^{\rm lip}$ the usual Lipschitz semi-norm and define $|f|_\g$ consequently.
\end{defi}

Note that $|\cdot|_{s}\le|\cdot
|_{s'}$ for $s\le s'$. Moreover the norms $|\cdot|_s$, $|\cdot|_{s,\g}$ satisfy the algebra,
interpolation and smoothing properties, namely for all $s\geq s_{0}$ 
there are $C(s)\geq C(s_{0})\geq1$ 
such that 
 if $A=A(\la)$ and $B=B(\la)$ depend on the parameter
 $\la\in\calI\subset\RRR$
in a Lipschitz way, then
\begin{subequations}
\begin{align}
|AB|_{s,\g}&\leq C(s)|A|_{s_0,\g}|B|_{s,\g}
+C(s_{0})|A|_{s,\g}|B|_{s_0,\g},\label{eq:2.11a}\\
|AB|_{s,\g}&\leq C(s)|A|_{s,\g}|B|_{s,\g}.\label{eq:2.11b}\\
||Ah||_{s,\g}&\leq C(s)(|A|_{s_0,\g}||h||_{s,\g}
+|A|_{s,\g}||h||_{s_0,\g}),\label{eq:2.13b}\\
|\Pi_{N}^{\perp}A|_{s,\g}&\leq N^{-\be}|A|_{s+\be,\g}, \quad \be\geq0,\label{eq:2.22}
\end{align}\label{stimenorma}
\end{subequations}
where
\begin{equation}\label{smoothop}
\left(\Pi_{N}A\right)_{k}^{k'}:=\left\{
\begin{aligned}
&A_{k}^{k'},\qquad \dist(k,k')\le N,\\
&0, \qquad \mbox{otherwise}.
\end{aligned}
\right.
\end{equation}
and $\Pi_N^{\perp}:=\uno-\Pi_N$.
The proof of the bounds \eqref{stimenorma} can be found in \cite{BB1} for the case
of the $s$-decay norm; given any norm $|\cdot|$ satisfying \eqref{stimenorma}
then also the corresponding $|\cdot|_\gamma$  satisfies \eqref{stimenorma}. 

\begin{rmk}
Note that, by \eqref{eq:2.13b} if a matrix $A$ has finite  norm $|A |_s$ then it is a bounded operator on $H^s$.
\end{rmk}

\begin{lemma}{(\cite{BP}-Lemma 7.1)}\label{moltiplicazione.matrici}
For any compact Lie group $G$ of dimension $\gotd$, 
consider $a,b,c\in H^{s}(\TTT^d\times G) $ with $a,b$ real valued. Then the multiplication operator with matrix
$$
B= \begin{pmatrix}
a(\f,x) & c(\f,x) \\ \bar c(\f,x)& b(\f,x)
\end{pmatrix}
$$
is self-adjont in $L^{2}$ and, for any $s>(d+\gotd)/2$,
$$
\|B_{\{i\}}^{\{i'\}}\|_0 \le C(s)\frac{\max(\|a\|_{s},\|b\|_{s},\|c\|_{s})}{\av{i-i'}^{s-(d+\gotd)/2}} \, , 
\quad \forall i, i' \in \ZZZ^d \times \Lambda_+ \, . 
$$
 In the case of $SU(2)$ we have $\gotd= 3$ and we deduce that
 $$
 |B|_s\leq C(s)\max(\|a\|_{s+\nu_0},\|b\|_{s+\nu_0}\|c\|_{s+\nu_0})\,,\quad \nu_0= (2d+5)/2\ .
 $$
Moreover for the central characters of $SU(2)$ the following multiplication rule holds:
\begin{equation}\label{multicaz}
\chi_h\chi_m=\sum_{k=0}^{\min(h,m)}\chi_{h+m-2k}\ .
\end{equation}
\end{lemma}

As explained in the introduction, the use of the off-diagonal decay norm is crucial in the proof of the
existence of solutions. For this reason, we find it convenient to use it also for the proof of stability
results; however, one could prove such stability results by simply using the operator norm: this only
requires a little more care in handling the small divisors.

We note (using also the regularity assumption on ${\mathtt f}$) that the operator $T$ defined in
\eqref{Ti} satisfies the following properties: 
\begin{subequations}
\begin{align}
& \text{\bf (T\"oplitz in time)} \ \qquad
T_{l,j,\se}^{l',j',\se'}= T_{j,\se}^{j',\se'}(l-l')
\label{hyp.topliz}\\
&  \text{\bf (Off-diagonal  decay)} \ \qquad 
 | T(w)|_{s-\nu_0}\le C(s)(1+ \|w\|_{s}) \, , 
\label{hyp.normaT} \\
&  \text{\bf (Lipschitz)} \ \qquad
 | T(w)-T(w')|_{s-\nu_0} \le C(s)(\|w-w'\|_s +\label{lip}\\
 &\qquad\qquad\qquad\qquad\qquad\qquad\qquad\qquad+(\|w\|_s + \|w'\|_s)\|w-w'\|_{s_0}) \, ,  
 \nonumber
\end{align}
for all $ \|w\|_{s_0}, \|w'\|_{s_0}\le 2$ and $ s_0+\nu_0 < s < q-2 $. 
\end{subequations}
\begin{rmk}\label{ceppissima}
For $ s_0+\nu_0 < s < q-2 $, \eqref{hyp.normaT} and \eqref{lip} imply
\begin{equation}\label{staceppa}
|T(w_\e)|_{s-\nu_0,\gamma,\calI}\leq C(s)(1+\|w_\e\|_{s,\gamma,\calI})\leq C(s)\ ;
\end{equation}
recall that $w_\e\in C^1(\calI,H^{s_1}(\gotK))$.
\end{rmk}

\section{The reduction algorithm}

It will be convenient to think of the equation \eqref{nls} as a Hamiltonian dynamical system on the phase
space $H^s(\calG)\times H^s(\calG)=H^s(\Lambda_+\times\gotA)$.

\begin{rmk}
Given a T\"oplitz in time matrix $T$ (see \eqref{hyp.topliz}), we can define, for all $\f\in \TTT^d$  a matrix on the phase space 
$H^s(\Lambda_+\times \gotA) $ by setting
 $$T_{j,\se}^{j',\gota'}(\f):=\sum_{l\in\ZZZ^d}T_{j,\se}^{j',\gota'}(l) e^{\ii l\cdot \f}$$
 and one has
 \begin{equation}\label{aaaaa}
\sup_{\f\in \TTT^d} |T(\f)|_{s}\leq C(s_{0})|T|_{s+s_{0}}.
\end{equation}
Note that in the l.h.s. we are considering the $s$-decay norm on $H^s(\Lambda_+\times\gotA)$ while in the r.h.s. we are considering the $s$-decay norm on $H^s(\gotK)$.
\end{rmk}

\begin{defi}[{\bf Hamiltonian vector field}]\label{hamvf}
Set $w=(u,\ol{u})\in H^s(\Lambda_+\times\gotA)$. We say that a vector field $X(w)$ is Hamiltonian
if there exists a real-on-real function $\calH(w)$ such that $X(w)=\ii J\nabla\calH(w)$ where
$$
J:=\begin{pmatrix}
0 & \uno \cr
-\uno & 0
\end{pmatrix}.
$$
In particular if $X(w)$ is linear, i.e. $X(w)=Mw$ for some matrix $M$, with
$$
M=\begin{pmatrix}
M^{+}_+ & M^+_- \cr
M^-_+ & M_-^-
\end{pmatrix},\qquad
M^+_+=\ol{M}^-_-=-(\ol{M}^+_+)^T\ \
M^+_-=(M^+_-)^T=\ol{M}^-_+
$$
then $X$ is Hamiltonian 
and the associated Hamiltonian function is
$$
\calH=-\frac{\ii}{2}\langle w, MJw\rangle
$$
and $\ii\s_3 M$ (see \eqref{pauli}) is a self-adjoint matrix.
\end{defi}

We now consider \eqref{vnls} linearised at $w_\e(\la)$ and we write it as a dynamical
system, namely we consider the linear equation
$$
\pDL_\e h =0
$$
where $\pDL_\e:=(\om\cdot\del_\f)+\ii \calL_\e$.
We want to show that
the linear operator $\pDL_\e$ can  be conjugated to a diagonal operator with purely imaginary spectrum.
We will also show that this change of variables acts on the phase space $H^s(\Lambda_+\times\gotA)$
by preserving the Hamiltonian structure.
Precisely we have the following result.

\begin{theorem}[{\bf KAM Theorem}]\label{KAMalgorithm}
Let ${\mathtt f}\in C^{q}$ for $q>s_0+\be+\nu_0+2$, $\g\in(0,1)$,
$s_2= \min(q-\be-\nu_0-2,s_1-\be-\nu_0)$
with $\be=6\tau+5$ for some $\tau>d$.
There exist constants $\epsilon_{0}$, $C$ such that,
if $\e\g^{-1}\leq\epsilon_{0}$,

%
then there exists a sequence of Lipschitz functions $\mu_j^\infty:\calI\to \RRR$
\begin{equation}\label{eq:4.4}
\mu_{j}^{\infty}(\la)=(j+\rho)^2-\rho^2+{\mathtt m}+ r_{j}^{\infty}(\la)\in\RRR, \quad \forall\;j\in\Lambda_+ ,\end{equation}
with
$|r_{j}^{\infty}|_{\g}\leq C\e$ for all $j\in\Lambda_+$ such that, setting
 \begin{equation}\label{eq:4.10}
 \begin{aligned}
\calS_{\infty}:=
\Big\{
\la\in\CCCC_\e &: |\la\tilde\om\cdot l\!+\!\gota\mu^{\infty}_{j}(\la)-\!\gota'\mu_{j'}^{\infty}(\la)|\geq
\frac{2\g}{\langle l\rangle^{\tau}},\\
& \;\forall l\in\ZZZ^{d}, \forall (j,\gota)\ne(j',\gota')\in\Lambda_+\times\gotA
\Big\}
\end{aligned}
\end{equation}
(where $\CCCC_\e$ is the set introduced in Theorem \ref{teoremone}),  the following holds.
For all $s_0\le s\le s_2$ and any $\la\in \calS_{\infty}$, there exists a bounded, invertible
linear operator $\Psi_{\infty}(\la) : { H}^{s}(\Lambda_+\times\gotA)\to { H}^{s}(\Lambda_+\times\gotA)$, with
bounded inverse $\Psi_{\infty}^{-1}(\la)$,
such that
\begin{equation}\label{eq:4.6}
\begin{aligned}
&\pDL_{\infty}(\la):=\Psi_{\infty}^{-1}(\la)\circ\pDL_\e\circ\Psi_{\infty}(\la)=
\la\tilde\om\cdot\del_{\f}\id+ \ii\DD_{\infty},\\
 &\DD_{\infty}:=\diag(\gota\mu_{j}^{\infty}(\la))_{\gota=\pm,j\in\Lambda_+}\,.
 \end{aligned}
\end{equation}
Moreover, the maps $\Psi_{\infty}(\la)$, $\Psi_{\infty}^{-1}(\la)$ satisfy
\begin{equation}\label{eq:4.8}
|\Psi_{\infty}(\la)-\id|_{s,\g,\calS_{\infty}}+
|\Psi_{\infty}^{-1}(\la)-\id|_{s,\g,\calS_{\infty}}\leq
\e \g^{-1} C(s)(1+||u_\e(\la)||_{s+\be+\nu_0,\g,\calS_{\infty}}).
\end{equation}
\end{theorem}


\subsection{The KAM step}\label{kamstep}

In this Section we  show in detail one step of the KAM iteration.

Let us consider a matrix  on the scale of spaces $H^s(\gotK)$
$$
\pDL= {\mathtt D}+  R
$$
with $\pD$ a diagonal matrix
$$
\begin{aligned}
&\pD=\diag (d_k)_{k\in\gotK}=
\diag(\ii\om\cdot l+\ii\gota\mu_j)_{\gota=\pm,j\in \Lambda_+,l\in \ZZZ^d}\,,\\
&\mu_j\in \RRR\,,\quad \mu_j=( j+\rho)^2-\rho^2
+{\mathtt m} +r_j  \,,\quad \sup_{j\in \Lambda_+}|r_j|< \infty
\end{aligned}
$$
and $\ii\sigma_3R$ is a self-adjoint, bounded matrix with finite $s$-decay norm for all $s< q-\nu_0$.
Moreover we assume that $R$ is T\"oplitz in time and, for all $\f\in \TTT^d$,
the vector field $ R(\f) w$
is Hamiltonian; see Definition \ref{hamvf}.

We construct a  canonical $\f$-dependent change of variables $\Phi$ which diagonalises $\pDL$
apart from  a {\em small remainder}; precisely $\Phi(\f)=e^{A(\f)}$ is the time-$1$ flow map generated by
a linear $\f$-dependent Hamiltonian system of the form
$$
\dot{x}= A(\f)x
$$
%
and we choose the matrix $A$ so that it solves the {\it homological} equation
\begin{equation}\label{homo}
\Pi_N R+ [A,\pD]=\text{ diag}(R)\,,
\end{equation}
where $[A,B]:=AB-BA$.
The smoothing operator $\Pi_{N}$ defined in \eqref{smoothop} is necessary for technical reasons: it is used
in order to obtain
suitable estimates on the high norms of the transformation $\Phi$, when the nonlinearity is merely
differentiable.

By \eqref{homo}, in the new variables we have the conjugated matrix
\begin{equation}\label{eq:4.1.22}
\begin{aligned}
& \pDL_1:=e^A\pDL e^{-A}= e^{{\rm ad}A}\pDL= 
\pD+\Pi_N R+ [A,\pD] + R_1= \pD_1+ R_1 \\
& \pD_1:=  \pD+\text{diag}(R) \\
&R_1:=
\Pi^\perp_N R+ \sum_{m\ge2}\frac{1}{m!}[A,\diag(R)-\Pi_N R]^{m-1}
+ \sum_{m\ge1}\frac{1}{m!}[A,R]^{m},
\end{aligned}
\end{equation}
where $[A,B]^m:=[A,[A,B]^{m-1}]$ and $\Pi_N^{\perp}:=\uno-\Pi_N$.
Note that $\pD_1$ has the same form as $\pD$; in particular
\begin{equation}\nonumber
\begin{aligned}
&\pD_1=\diag (d_k)_{k\in\gotK}=
\diag(\ii\om\cdot l+\ii\gota\mu_j^{(1)})_{\gota=\pm,j\in \Lambda_+,l\in \ZZZ^d}\,,\\
&\mu_j^{(1)} \in \RRR\,,\quad \mu_j^{(1)}=( j+\rho)^2-\rho^2+{\mathtt m} +r^{(1)}_j  \,,
\quad \sup_{j\in \Lambda_+}|r_j^{(1)}|< \infty.
\end{aligned}
\end{equation}

In order  to solve the {homological} equation \eqref{homo} we simply note that
\begin{equation}\label{eq:2.4}
[A,\pD]_{k}^{k'}=A_{k}^{k'} (d_{k'}-d_{k})
\end{equation}
and hence we can set
\begin{equation}\label{sol.homo}
A_{k}^{k'}=\left\{
\begin{aligned}
&\frac{R_{k}^{k'}}{d_{k}-d_{k'}}, \qquad 0<\dist(k,k')\le N,\\
&0, \qquad \mbox{ otherwise}.
\end{aligned}
\right.
\end{equation}
Moreover, defining
$$
\calS_+:=\{\la\in\calS\;:\; |d_{l,j,\gota}-d_{l',j',\gota'}|> \g |l-l'|^{-\tau}\,, \mbox{ for all } 0<|l-l'|\le N \}
$$
one has the bound
\begin{equation}\label{eq:4.1.33}
|A|_{s,\g,\calS_+}\leq C N^{2\tau+1}\g^{-1}
|R|_{s,\g,\calS}\,.
\end{equation}
Finally $A$ is T\"oplitz in time and $\Phi$ is a canonical change of variables provided that
$| R|_{s}$ is small enough for some $s$.

The eigenvalues $\mu_{j}^{(1)}$ satisfiy
\begin{equation}\label{45}
|\mu_{j}^{(1)}-\mu_{j}|^{{\rm lip}}=|r_{j}^{(1)}-r_{j}|^{{\rm lip}}=|\diag(R)|^{{\rm lip}}\leq
|R|^{{\rm lip}}_{s_{0}}, \quad j\in\Lambda_+,
\end{equation}
while the remainder $R_{1}$ satisfies
\begin{equation}\label{46}
\begin{aligned}
&|R_1|_{s,\g}\leq C({s})( N^{-\be}|R|_{s+\be,\g}+N^{2\tau+1}\g^{-1}
|R|_{s,\g}|R|_{s_{0},\g}),\\
&|R_1|_{s+\be,\g}\leq C({s+\be})( |R|_{s+\be,\g}+N^{2\tau+1}\g^{-1}
|R|_{s+\be,\g}|R|_{s_{0},\g}).
\end{aligned}
\end{equation}

\subsection{The iterative Lemma}

We now iterate the procedure above infinitely many times. Throughout the procedure we shall keep track
of the parameter $\g$ since eventually we want to fix it so that it is small with $\e$. On the other hand
we will systematically ignore the constants not depending on the iteration step, $\e$ and $\g$.

\begin{lemma}\label{teo:KAM}
Let $q>s_{0}+\be+\nu_0+2$ and set $\calS_{0}:=\CCCC_\e$ and
$\pDL_0:=(\om\cdot\del_\f)+\ii\calL_\e$; see \eqref{line}. There exist a constant $C_{0}>0$ and $N_{0}\in\NNN$ large
(independent of $\e,\g$), such that
if
\begin{equation}\label{eq:4.15}
\e N_{0}^{C_{0}}\g^{-1}|T|_{s_{0}+\be,\g}\leq1,
\end{equation}
then, for any $n\geq1$, if we set $N_n:=N_0^{(\frac{3}{2})^n}$ the following holds.

\noindent
$({\bf S1})_{n}$ 
Setting
\begin{equation}\label{eq:419bis}
\begin{aligned}
\calS_{n}:=
\Big\{\la\in\calS_{n-1} &: 
|\la\tilde\om\cdot l\!+\! \gota\mu_{j}^{(n-1)}(\la)\!-\! \gota'\mu_{j'}^{(n-1)}(\la)|\geq
\frac{\g}{\langle l\rangle^{\tau}},\\
&\quad\ \forall\, | l|\leq N_{n-1}, \, (j,\se)\ne(j',\se')\in\Lambda_+\times\gotA
\Big\},
\end{aligned}
\end{equation}
then, for all $\la\in\calS_n$, we can apply the KAM step described in Section \ref{kamstep}
to $\pDL_{n-1}$, namely there exists a T\"oplitz in time matrix $A_{n-1}$ which defines a
canonical change of variables $\Phi_{n-1}:H^s(\Lambda_+\times\gotA)\to H^s(\Lambda_+\times\gotA)$ with
$\Phi_{n-1}:={e}^{A_{n-1}}$ such that
\begin{equation}\label{eq:4.16}
\begin{aligned}
\pDL_{n}&:=\Phi_{n-1}^{-1} \pDL_{n-1}\Phi_{n-1}:=\om\cdot\del_{\f}  +\ii\DD_{n}+R_{n},\\
\qquad
\DD_{n}&={\rm diag}(\gota\mu_j^{(n)})_{\gota=\pm,j\in\Lambda_+},
 \\
\mu_j^{(n)}&=\mu_{j}^{(n)}(\la)=(j+\rho)^2-\rho^2+{\mathtt m}+r_{j}^{(n)}(\la)\in\RRR,
\end{aligned}
\end{equation}
with
\begin{equation}\label{eq:4.20}
|r_{j}^{(n)}|_{\g}:=|r_{j}^{(n)}|_{\g,\calS_{n}}\leq\e C,
\end{equation}
and the vector field $R_{n}$ is Hamiltonian.

\noindent
$({\bf S2})_{n}$
The matrix $A_{n-1}$  satisfies
\begin{equation}\label{eq:4.22}
\begin{aligned}
|A_{n-1}|_{s,\g}\leq 
\e|T|_{s+\be,\g} N_{n-1}^{2\tau+1}N_{n-2}^{-\be+1}.
\end{aligned}
\end{equation}

\noindent
$({\bf S3})_{n}$
For all $s\in[s_{0},s_2]$ one has
\begin{equation}\label{eq:4.21}
\begin{aligned}
|R_n|_{s,\g}&\leq \e|T|_{s+\be,\g}N_{n-1}^{-\be+1},\\
|R_n|_{s+\be,\g}&\leq \e|T|_{s+\be,\g}N_{n-1},
\end{aligned}
\end{equation}

\noindent
$({\bf S4})_{n}$ For all $j\in\NNN$ there exists  Lipschitz extensions
$\tilde{\mu}_{j}^{(n)}(\cdot) : \calI \to \RRR$ of $\mu_{j}^{(n)}(\cdot):\calS_{n}\to \RRR$,
such that one has
\begin{equation}\label{eq:4.23}
|\tilde{\mu}_{j}^{(n)}-\tilde{\mu}_{j}^{(n-1)}|_{\g}\leq
|R_{n-1}|_{s_{0},\g}.
\end{equation}

\end{lemma}

\noindent
{\it Sketch of the proof.}
We proceed by induction. The case $n=1$ follows by the smallness hypothesis.
Indeed \eqref{eq:4.15} implies the smallness of $|R_0|_{s_0+\beta,\g}$ which in turn
by \eqref{eq:4.1.33} implies that $|A_{0}|_{s_0+\be,\g}<1/2$. Then $\Psi_0$ is well defined
and the bounds \eqref{eq:4.20}, \eqref{eq:4.21} and \eqref{eq:4.22}
as well as $({\bf S4})_{1}$ follow by \eqref{46} and \eqref{45}. Recall that by the Kirszbraun Theorem we
can extend $r^{(1)}_j(\la)$ to a Lipschitz function on the whole interval $\calI=[1/2,3/2]$.

For $n\ge2$ we start by defining $A_{n-1}$ on the set $\calS_{n}$ using the homological
equation \eqref{homo} with $A\rightsquigarrow A_{n-1}$ and $R\rightsquigarrow R_{n-1}$; in particular
$r_j^{(n)}=r_{j}^{(n-1)}+(\diag(R_{n-1}))_j$ so that the bound \eqref{eq:4.20} follows by the inductive hypothesis
and hence $({\bf S1})_{n}$ follows. Then \eqref{eq:4.1.33}
together with $({\bf S3})_{n-1}$ directly implies $({\bf S2})_{n}$. To prove $({\bf S3})_{n}$ we use \eqref{46}
and $({\bf S3})_{n-1}$; precisely we have that $R_n$ is defined as in \eqref{eq:4.1.22}
(with clearly $R_1\rightsquigarrow R_n$) and hence it satisfies the bound \eqref{46} with $N\rightsquigarrow N_{n-1}$.
But then we may use $({\bf S3})_{n-1}$ and obtain
$$
\begin{aligned}
|R_n|_{s,\g}&\le C(s) \e |T|_{s,\g}(N_{n-1}^{-\be} N_{n-2}+\g^{-1}\e |T|_{s_0+\be,\g}N_{n-1}^{2\tau+1}
N_{n-2}^{-2\be+2})\\
&\stackrel{\eqref{eq:4.15}}{\le} C(s)\e |T|_{s,\g}(N_{n-2}^{-\frac{3}{2}\be+1}+N_{n-2}^{3\tau-2\be+\frac{7}{2}})
\end{aligned}
$$
which implies the desired bound since $\be=6\tau+5$. The second bound in \eqref{eq:4.21} follows similarly.
Finally $({\bf S4})_{n}$ follows by the Kirszbraun Theorem.
\EP

\subsection{Proof of Theorem \ref{KAMalgorithm}}

First we verify that the hypotheses of Theorem \ref{KAMalgorithm} imply those of Lemma \ref{teo:KAM}.
Indeed, since $ s_0+\nu_0 < s_0+\be < q-2 $, we can apply Remark \ref{ceppissima}. Then, by \eqref{staceppa} we have:
$$
\e N_{0}^{C_{0}}\g^{-1}|T|_{s_{0}+\be,\g}\leq \e C (s_0+\be)N_{0}^{C_{0}}\g^{-1} ;
$$
recalling that $N_0,C_0,C(s_0+\be)$ are independent of $\e,\g$ this is smaller than $1$ provided $\e\g^{-1}$ is small enough, which amounts to taking $\epsilon_0$ small in Theorem \ref{KAMalgorithm}.

Now we have to prove that the iteration described in Lemma \ref{teo:KAM} converges.
We show that there exists the ``limit'' change of variables $\Psi_{\infty}$. 
For any $\la\in\cap_{n\geq0}\calS_{n}$ we define
\begin{equation}\label{eq:4.1.3}
\Psi_{n}:=\Phi_{0}\circ\Phi_{1}\circ\ldots\circ\Phi_{n}
\end{equation}
and we note that $\Psi_{n+1}=\Psi_{n}\circ\Phi_{n+1}$.
Then, one has
\begin{equation}\label{eq:4.1.4}
|\Psi_{n+1}|_{s_{0},\g}\stackrel{(\ref{eq:2.11b})}{\leq}
|\Psi_{n}|_{s_{0},\g}+C|\Psi_{n}|_{s_{0},\g}|\uno-\Phi_{n+1}|_{s_{0},\g}
\stackrel{(\ref{eq:4.22})}{\leq} 
|\Psi_{n}|_{s_{0},\g}(1+\de^{(0)}_{n}),
\end{equation}
where we used that $\Phi_n=e^{A_n}$ and we have defined
\begin{equation}\label{eq:4.1.5}
\de_{n}^{(s)}:=\e K \g^{-1}N_{n+1}^{2\tau+1}N_{n}^{-\be+1}|T|_{s,\g}, 
\end{equation}
for some constant $K>0$. Now, by iterating (\ref{eq:4.1.4}) and using the (\ref{eq:4.15}), (\ref{eq:4.22}),
we obtain
\begin{equation}\label{eq:4.1.6}
|\Psi_{n+1}|_{s_{0},\g}\leq|\Psi_{0}|_{s_{0},\g}\prod_{n\geq0}(1+\de_{n}^{(s_0)})
\leq2\,.
\end{equation}
The estimate on the high norm follows by interpolation and one obtains
\begin{equation}\label{eq:4.1.7}
|\Psi_{n+1}|_{s,\g}{\leq}C(s)\left(1+\e\g^{-1}|T|_{s+\be,\g}\right).
\end{equation}
Thanks to (\ref{eq:4.1.7}) one easily sees that 
 the sequence $\Psi_{n}$ is a Cauchy sequence w.r.t. the norm $|\cdot|_{s,\g}$; in particular one
 has
\begin{equation}\label{eq:4.1.8}
|\Psi_{n+m}-\Psi_{n}|_{s,\g}\leq C(s)
\e\g^{-1}|T|_{s+\be,\g}N_{n}^{-1}.
\end{equation}
As a consequence one has that $\Psi_{n}\stackrel{|\cdot|_{s,\g}}{\to}\Psi_{\infty}$ and
 (\ref{eq:4.8}) is verified.

Let us now define for $j\in\Lambda_+$,
\begin{equation}\label{eq:4.1.1}
\mu_{j}^{\infty}:=\lim_{n\to+\infty}\tilde{\mu}^{(n)}_{j}(\la)
=(j+\rho)^2-\rho^2+{\mathtt m}+
\lim_{n\to+\infty}\tilde{r}_{j}^{(n)}(\la)
\end{equation}
and note that, for any $n\in \NNN$, $j\in\Lambda_+$, one has
\begin{equation}\label{eq:4.1.10}
\begin{aligned}
|\mu_{j}^{\infty}-\tilde{\mu}_{j}^{(n)}|_{\g,\calI}&\leq
\sum_{m=n}^{\infty}|\tilde{\mu}^{(m+1)}_{j}-\tilde{\mu}_{j}^{(m)}|_{\g,\calI}
\stackrel{\eqref{eq:4.15},(\ref{eq:4.23}),(\ref{eq:4.21})}{\leq}
\g N_{n-1}^{-\be+1}.
\end{aligned}
\end{equation}
Hence we have proved that for all $\la\in\cap_{n\ge 0}\calS_n$ the linear operator
$\pDL_\e$ is conjugated via $\Psi_\io$ to $\pDL_\io$; see \eqref{eq:4.6}.

In order to conclude the proof of Theorem \ref{KAMalgorithm}
we only need to prove that
\begin{equation}\label{eq:4.1.12}
\calS_{\infty}\subseteq\bigcap_{n\geq0}\calS_{n}.
\end{equation}
We show by induction that for any $n>0$ then 
$\calS_{\infty}\subseteq \calS_{n}$.
By definition we have $\calS_{\infty}\subseteq\calS_{0}:=\CCCC_\e$.
Assume that $\calS_{\infty}\subseteq \cap_{p=0}^{n} \calS_{p}$, so that the
$\mu^{(n)}_{j}$'s are well defined and coincide with their extension.
Then, for any fixed  $(j,\gota),(j',\se')\in\Lambda_+\times\gotA$, and any $l\in\ZZZ^d$
 we have
\begin{equation}\label{eq:4.1.13}
|\om\cdot l+\gota\mu_{j}^{(n)}-\se'\mu_{j'}^{(n)}|
\stackrel{(\ref{eq:4.10}),(\ref{eq:4.1.10})}{\geq}\frac{2\g}{\langle l\rangle^{\tau}}
-2\g N_{n-1}^{-\be+1}.
\end{equation}

Now, since  $|l|\leq N_{n}$ and $\be=6\tau+5$, we have
\begin{equation}\label{eq:4.1.15}
|\om\cdot l+\se\mu_{ j}^{(n)}-\se'\mu_{j'}^{(n)}|\geq\frac{\g}{\langle l\rangle^{\tau}},
\end{equation}
which implies $\calS_\io\subseteq \cap_{p=0}^{n+1} \calS_{p}$. Hence the assertion follows.
\EP

\subsection{Measure estimates}\label{puffetta}

We define the set of ``resonant parameters'', namely
\begin{equation}\label{reso}
\begin{aligned}
\calmR&:=\bigcup_{l\in\ZZZ^d}\bigcup_{\substack{j,j'\in\Lambda_+\\ \se,\se'\in\gotA\\ (j,\se)\ne(j',\se')}}
\calmR_{l,j,j',\gota,\se'}\,,\\
&\calmR_{l,j,j',\gota,\se'}:=
\left\{
\la\in\calI : |\la\tilde\om\cdot l\!+\!\gota\mu^{\infty}_{j}(\la)-\!\gota'\mu_{j'}^{\infty}(\la)|\leq
\frac{2\g}{\langle l\rangle^{\tau}}
\right\}
\end{aligned}
\end{equation}
and we want to prove that $\meas(\calmR)=O(\g)$; clearly this implies that
$\meas(\calS_\io)\ge\meas(\CCCC_\e)-O(\g)=1-C_0\e^{1/S} -C_1\g$ so that by choosing
$\g = \e^{1/S}$ one has $\meas(\calS_\io)\to1$ as $\e\to0$.
Note that this choice of $\g$ is compatible with the smallness condition \eqref{eq:4.15}.

First of all we note that $(j,\se)\ne(j',\se')$ implies
$$
|\gota\mu^{\infty}_{j}(\la)-\!\gota'\mu_{j'}^{\infty}(\la)| \geq \frac58 -C\e\ .
$$
Then, for all $(j,\se)\ne(j',\se')$ the condition
$$
|\la\tilde\om\cdot l\!+\!\gota\mu^{\infty}_{j}(\la)-\!\gota'\mu_{j'}^{\infty}(\la)|\leq
\frac{2\g}{\langle l\rangle^{\tau}}
$$
implies that
$$
|\tilde\om\cdot l|\geq \frac23 \left(\frac58 -C\e -2\g\right) \ge\frac 13\ .
$$
This means that if $|\tilde\om\cdot l|<1/3$, then $\calmR_{l,j,j',\gota,\se'}=\emptyset$. Otherwise if $|\tilde\om\cdot l|\geq1/3$, one has (since the $\mu_j^{\io}$'s are
Lipschitz functions on the whole interval $\calI$)
$$
|\tilde\om\cdot l|-2\sup_{j\in\Lambda_+}|\mu^{\infty}_{j}(\la)|^{\rm lip}\ge 
\frac13-C\e\g^{-1}\ge\frac14\ ,
$$
which implies the measure estimate
$$
\meas(\calmR_{l,j,j',\se,\se'})
\le 8\g\langle l\rangle^{-\tau}\ .
$$

Now we claim that
\begin{equation}\label{fanculo}
|\se(j+\rho)^2-\se'(j'+\rho)^2|>6|l|
\end{equation}
implies $\calmR_{l,j,j',\se,\se'}=\emptyset$. For $l=0$ this is trivial by the definition of $\mu_j^\io$.
For $l\ne0$ \eqref{fanculo} implies
$$
|\se\mu_j^{\io}-\se'\mu_{j'}^{\io}|>3|l|\ge 2|\la||\tilde{\om}|_1|l|
$$
and our claim follows.

Finally, the negation of \eqref{fanculo} implies $|j|,|j'|<9|l|$ so that we can bound
$$
\meas(\calmR)\le C\g\sum_{l\in\ZZZ^d} \langle l\rangle^{2-\tau}
$$
and the wanted measure estimate follows for $\tau>d+2$.

\section{Proof of Theorem \ref{thm.reduci}}
Assume that $q>(15d+57)/2$ and take $s_0=1+d/2$ in Definition \ref{def:Ms}. Then set $\tau=d+3$ so that $\alpha= 7d+26$, $s_2=\min(q-7d-55/2,s_1-7d-26)$. It is easily seen that these choices of parameters satisfy all the constraints in Theorem \ref{KAMalgorithm}. Fix $\gamma=\e^{1/S}$ with $S$ given in Theorem \ref{teoremone}. Since with this choice $\e\g^{-1}$ is small with $\e$, then Theorem \ref{thm.reduci} follows from Theorem \ref{KAMalgorithm} by choosing $\DD=\DD_\infty$, $\Psi=\Psi_\infty$ and $\mathcal S=\mathcal S_\infty$. The measure estimate \eqref{pillola} follows by Section \ref{puffetta} since the complementary to $\mathcal S_\infty$ in $\cal I$ is $\calmR$.

Finally, in order to prove \eqref{scaramouche}, we
 observe that
\begin{eqnarray}\label{fandango}
& & \Big|\|h(t)\|_{s_2}-\|h(0)\|_{s_2}\Big|\leq\\
\nonumber & & \leq\Big|\|h(t)\|_{s_2}-\|v(t)\|_{s_2}\Big|+\Big|\|v(t)\|_{s_2}-\|v(0)\|_{s_2}\Big|+\Big|\|v(0)\|_{s_2}-\|h(0)\|_{s_2}\Big|\ .
\end{eqnarray}
Now, the second term in the r.h.s. of \eqref{fandango} is identically zero, while the first and the third can be estimated via \eqref{poropo}, obtaining
$$
\Big|\|h(t)\|_{s_2}-\|h(0)\|_{s_2}\Big|\leq C \e^a (1 +\|u_\e(\la)\|_{s_2+\alpha})(\|h(t)\|_{s_2}+\|h(0)\|_{s_2})
$$
which implies \eqref{scaramouche}.

\section{Final remarks and open problems}\label{salaminchia}

For the sake of simplicity, we confined ourselves to the case of $SU(2)$ and $SO(3)$. However, the only important conditions are the fact that $\Lambda_+$ is one-dimensional and that the eigenspaces of the Laplacian restricted to the subspace of central function are one-dimensional. This means that our results extends directly to the case of spherical varieties of rank $1$ provided that we restrict ourselves to symmetric functions.

It seems extremely reasonable that most results that hold true for tori can be extended to the case of homogeneous manifold (indeed the harmonic analysis is very similar), provided that one restricts him/herself to central functions in order to avoid multiplicity of the eigenvalues.

A very natural question is whether the reducibility results by Eliasson-Kuksin \cite{EK1,EK} (at least in the simplified case considered in Procesi-Xu \cite{PX}) can be extended also to this setting. In other words, this would mean to be able to extend the result of the present paper to the case of arbitrary rank. Of course, the KAM scheme works regardless of the rank: the problem concerns only the measure estimates. Indeed, in the case of rank greater than $1$, equation \eqref{fanculo} does not imply $|j|,|j'|<9|l|$ and hence the union in \eqref{reso} may cover the whole interval $\calI$. In order to overcome this difficulty, one needs more precise information on the eigenvalue asymptotics. This would require a suitable extension of the notion of T\"oplitz-Lipschitz or quasi-T\"oplitz matrices, which is most probably feasible but technically extremely complicated.

Another interesting problem would be to consider also autonomous equations. This is related to a better understanding of the Birkhoff normal form on compact manifolds. This is still an open problem, except for the case of tori and Zoll manifolds. Naturally, there should be no problem in the case of $SU(2)$ or objects of rank $1$. For more general Lie groups, in principle one can compute the Birkhoff normal form by using the eigenfunction multiplication rules (\cite{BP}, formula (2.20)); however, it would probably require some very heavy computations and it is not clear which kind of information one can obtain in this way.

This research was supported by the European Research Council under
FP7 ``Hamiltonian PDEs and small divisor problems: a dynamical systems approach''


\begin{thebibliography}{99} 
 
{\small 


\bibitem{BBM1}{
Baldi P., Berti M., Montalto R.,
\textit{KAM for quasi-linear and fully non-linear forced perturbations of Airy equation},
Math. Annalen 359 (2014), 471-536
 }


\bibitem{BBM3}{
Baldi P., Berti M., Montalto R.,
\textit{KAM for quasi-linear KdV},
C.R. Math. Acad. Sci. Paris 352 (2014), no.7-8, 603-607
 }


\bibitem{BBM2}{
Baldi P., Berti M., Montalto R.,
\textit{KAM for autonomous quasi-linear perturbations of KdV},
 preprint, 2014
 }


\bibitem{BB1}{ 
Berti M., Bolle Ph., 
\textit{Quasi-periodic solutions with Sobolev regularity of
NLS on $\TTT^{d}$ with a multiplicative potential},  Journal European Math. Society, 
15, 229--286, 2013.} 
 
\bibitem{BB2}{ 
Berti M., Bolle Ph., 
\textit{Sobolev quasi periodic solutions of multidimensional
wave equations with a multiplicative potential},
Nonlinearity, 25, 2579--2613, 2012.} 
 
\bibitem{BBP}{ 
Berti M., Bolle Ph., Procesi M.
\textit{An abstract Nash-Moser theorem with parameters and applications to PDEs},
Ann. I. H. Poncar\'e, 1, 377--399, 2010.} 

\bibitem{BCP}{
Berti M., Corsi L., Procesi M.,
\textit{An abstract Nash-Moser theorem
and quasi-periodic solutions for NLW and NLS on 
compact Lie groups and homogeneous manifolds},
 Comm. Math. Phys. published online Aug. 2014 } 

\bibitem{BP}{
Berti M., Procesi M.,
\textit{Nonlinear wave and Schr\"odinger equations on compact Lie groups and
Homogeneous spaces}, Duke Math. J., 159, 479--538, 2011.
}


\bibitem{B1}  Bourgain J., {\it  Construction of quasi-periodic solutions
for Hamiltonian perturbations of linear equations and applications
to nonlinear PDE}, Internat. Math. Res. Notices, no. 11, 1994.


\bibitem{B3}  Bourgain J., {\it Quasi-periodic solutions of Hamiltonian
perturbations of $2D$ linear Schr\"odinger equations},
Annals of Math. 148, 363-439, 1998.

\bibitem{B5}  Bourgain J.,  {\it Green's function estimates for lattice Schr\"odinger
operators and applications}, Annals of Mathematics Studies 158,
Princeton University Press, Princeton, 2005.


\bibitem{CY} Chierchia L., You J., {\it KAM tori for 1D
nonlinear wave equations with periodic boundary conditions}, Comm.
Math. Phys. 211, 497-525, 2000.


\bibitem{CW} Craig W., Wayne C. E., {\it Newton's method and periodic solutions
of nonlinear wave equation}, Comm. Pure  Appl. Math. 46,
1409-1498, 1993.


\bibitem{EK1} Eliasson L. H., Kuksin S., {\it On reducibility of Schr\"odinger equations with
quasiperiodic in time potentials}, Comm. Math. Phys, 286, 125-135, 2009.

\bibitem{EK} Eliasson L.H., Kuksin S.,
{\it KAM for non-linear Schr\"odinger equation},  Annals of Math.,
172, 371-435, 2010.

\bibitem{FP}
Feola R., Procesi M.
{\it Quasi-periodic solutions for fully non-linear forced NLS},
preprint, 2014


\bibitem{FS} Fr\"ohlich, J., Spencer, T.
{\it Absence of diffusion in the Anderson tight binding model for large disorder or low energy}, 
Comm. Math. Phys. 88 (1983), no. 2, 151--184. 


\bibitem{GXY} Geng J., Xu X.,  You J., {\it An infinite dimensional KAM theorem and 
its application to the two dimensional cubic 
Schr\"odinger equation}, Adv. Math. 226, 5361-5402, 2011.


%


\bibitem{K1} Kuksin S., {\it Hamiltonian perturbations of infinite-dimensional linear systems with imaginary spectrum},
Funktsional Anal. i Prilozhen., 21, 22-37, 95, 1987.

\bibitem{K2} Kuksin S., {\it Analysis of Hamiltonian PDEs}, Oxford
Lecture series in Mathematics and its applications 19, Oxford University
Press, 2000.

\bibitem{KP} Kuksin S., P\"oschel J.,
{\it Invariant Cantor manifolds of quasi-periodic oscillations for a nonlinear Schr\"odinger
equation},  Ann. Math. (2) 143 149-79, 1996.

\bibitem{M}{
J. Moser,
\emph{Convergent series expansions for quasi--periodic motions},
Math. Ann.
\textbf{169} (1967) 136--176. }

\bibitem{Po2}{ P\"oschel J.,
{\it A KAM-Theorem for some nonlinear PDEs},  Ann. Sc. Norm. Pisa, 
23, 119-148, 1996.}



\bibitem{Po3}{ P\"oschel J.,
{\it On elliptic lower-dimensional tori in Hamiltonian systems.},   
Math. Z. 202 (1989), no. 4, 559--608.}

\bibitem{PP} {Procesi M., Procesi C., {\it A KAM algorithm for the resonant nonlinear
Schr\"odinger equation}, preprint 2013.}

\bibitem{PX} {Procesi M., Xu X., {\it Quasi-T\"oplitz Functions in KAM Theorem}, 
 SIAM J.Math. Anal. 45, 4, 2148-2181, 2013.}
 
\bibitem{W1}  {Wayne E., {\it Periodic and quasi-periodic solutions
of nonlinear wave equations via KAM theory}, Comm. Math. Phys.
127, 479-528, 1990.}







}
\end{thebibliography}
\end{document}